\documentclass[11pt]{article}
\usepackage[T1]{fontenc}
\usepackage{lmodern}
\usepackage{amsmath,amssymb,amsthm}
\usepackage[margin=1in]{geometry}
\usepackage{microtype}
\usepackage{needspace}
\usepackage[bottom]{footmisc}
\usepackage{xurl}
\usepackage{xcolor}
\definecolor{footnoteblue}{RGB}{0,70,160}
\usepackage[hidelinks]{hyperref}
\makeatletter
\renewcommand{\@makefnmark}{\hbox{\@textsuperscript{\normalfont\color{footnoteblue}\@thefnmark}}}
\makeatother
\hypersetup{pdftitle={Singular Loci of SL2(C)-Representation Varieties of Non-Hopfian Baumslag--Solitar Groups},
 pdfauthor={S. Liriano},
 pdfsubject={Representation varieties, singular loci, and topology}}
\usepackage{booktabs}
\usepackage{fancyhdr}
\newtheoremstyle{compactplain}{6pt}{6pt}{\itshape}{}{\bfseries}{.}{0.5em}{}
\theoremstyle{compactplain}
\newtheorem{theorem}{Theorem}[section]
\newtheorem{lemma}[theorem]{Lemma}
\newtheorem{proposition}[theorem]{Proposition}
\newtheorem{corollary}[theorem]{Corollary}
\newtheorem*{paritylemma}{Parity lemma}
\DeclareMathOperator{\SL}{SL}
\DeclareMathOperator{\GL}{GL}
\DeclareMathOperator{\Hom}{Hom}
\DeclareMathOperator{\Sing}{Sing}
\DeclareMathOperator{\tr}{tr}
\newcommand{\CC}{\mathbb C}
\newcommand{\FF}{\mathbb F}
\newcommand{\ZZ}{\mathbb Z}
\newcommand{\Aff}{\mathbb A}
\newcommand{\RR}{R_{m,n}}
\newcommand{\Singmn}{S_{m,n}}

\title{\Large Singular Loci of $\SL_2(\CC)$-Representation Varieties\\[0.25em]
of Non-Hopfian Baumslag--Solitar Groups}
\author{S.~Liriano}
\date{November 2, 2026}

\begin{document}
\maketitle

\begin{abstract}
We study the singular loci of the reduced $\SL_2(\CC)$-representation
varieties $R_{m,n}$ of non-Hopfian Baumslag--Solitar groups
$B(m,n)=\langle a,b\mid b^{-1}a^m b=a^n\rangle$.
We determine their dimensions and profile functions, which list the
numbers of irreducible components in descending order of dimension.
The components have dimensions three or five.
Their reduced singular loci have the explicit description
\[
 \Sing R_{m,n}
 =\left\{(\varepsilon I,Y):
   \varepsilon^m=\varepsilon^n,\quad
   (\tr Y)^2=\frac{(m+n)^2}{mn}\right\},
\]
where $\varepsilon\in\{\pm1\}$ and $Y\in\SL_2(\CC)$.
Each is a disjoint union of $\mu$ copies of the smooth affine quadric
$Q=\{(x,y,z)\in\CC^3:x^2+yz=1\}$, where $\mu=2$ for $m-n$ odd
and $\mu=4$ for $m-n$ even. Thus its complex algebraic isomorphism
type depends only on this parity.
Each quadric is homotopy equivalent to $S^2$, so the singular locus
has integral homology $\ZZ^\mu$ in degrees zero and two and zero
in all other degrees. An algebraic stratification gives its compactly
supported E-polynomial $\mu(q^2+q)$, with $q=uv$.
At every singular point, $R_{m,n}$ is locally analytically the product
of a smooth surface and a nodal curve.
\end{abstract}

\section{Introduction}

For a fixed complex affine algebraic group $A$, the representation variety
$R_A(G)=\Hom(G,A)$, considered up to algebraic isomorphism, is an algebraic
invariant of the finitely generated group $G$, independent of the chosen
finite generating set.\footnote{See Liriano--Majewicz \cite[p.~1]{LM}.
An isomorphism $\phi:G\to H$ induces an algebraic isomorphism
$\phi^*:R_A(H)\to R_A(G)$, $\rho\mapsto\rho\circ\phi$, with inverse
$(\phi^{-1})^*$; in generator coordinates, both maps are given by regular
word maps.} Its irreducible components,
singularities, and topology provide ways to examine $G$ via algebraic geometry. In this paper we consider
\[
 B(m,n)=\langle a,b\mid b^{-1}a^m b=a^n\rangle,
 \qquad
 \RR=\Hom(B(m,n),\SL_2(\CC)),
\]
and restrict throughout to the non-Hopfian Baumslag--Solitar groups~\cite{BS62}.
Equivalently, $|m|,|n|>1$ and $m,n$ have different sets of prime divisors.\footnote{For this criterion, see Levitt \cite[p.~3]{Levitt}. In particular, the non-Hopfian hypothesis implies $|m|\ne|n|$.}
We give $\RR$ its reduced affine variety structure.
With the convention $a\mapsto X$, $b\mapsto Y$, its points are
ordered pairs $(X,Y)\in\SL_2(\CC)^2$ satisfying $X^mY=YX^n$.
Note, we are not taking the quotient by conjugation.
Our principal goal is the reduced singular locus $\Singmn=\Sing\RR$.

All homotopy types and singular homology groups in this paper refer to the
usual complex analytic topology. The algebraic and analytic local rings at
a complex point have canonically isomorphic completions. Consequently,
local dimension and nonsingularity agree under analytification.\footnote{Serre \cite[Section~2, no.~6, Proposition~3 and Corollary~2]{Serre}.}

The starting point for this work was a series of unpublished computations
involving Gr\"obner bases in representation varieties over finite fields
conducted over several years, via SageMath~\cite{Sage},
Singular~\cite{Singular} (included in SageMath), and
Macaulay2~\cite{M2}. Those experiments included counts of rational points on the
representation varieties and on their singular loci. They led us to suspect
that the singular loci in the present family should have dimension two,
even when the dimensions of the ambient representation varieties varied.
The arithmetic calculations supplied the intuition for the geometric
statements proved here. In earlier unpublished work of more than two decades ago, the author had
obtained dimension calculations and the reducibility status of these
representation varieties. Nevertheless, in the sequel we make use of
existing literature to arrive at many of our conclusions.

There is a classical precedent for this passage from arithmetic evidence
to topological questions. In his 1949 paper on equations over finite
fields, Andr\'e Weil proposed a connection between point counts over finite
extensions and the Betti numbers of corresponding smooth complex projective
varieties.\footnote{Weil \cite[pp.~498, 507--508]{Weil}. This is a paraphrase of the guiding principle, not a quotation. The relevant topological comparison concerns complex Betti numbers and excludes exceptional primes.}
Katz's polynomial-count theorem also applies to affine or singular
varieties: if a model $\mathcal V$ of a complex variety $V$ over a
finitely generated subring $\mathcal A\subset\CC$ has
$\#\mathcal V_\varphi(\FF_q)=C_V(q)$ for every homomorphism
$\varphi:\mathcal A\to\FF_q$, with one polynomial
$C_V(T)\in\ZZ[T]$, then $E(V;u,v)=C_V(uv)$.\footnote{Katz,
Appendix to Hausel--Rodriguez-Villegas \cite[Theorem~6.1.2(3),
pp.~618--619]{HRV}.}
For our singular loci, the quadric decomposition below supplies such
a model after adjoining $\sqrt{mn}$ and inverting $2mn(m-n)(m+n)$.
The finite-field calculations thus provide questions and potential
evidence, together with a precise arithmetic interpretation of the
E-polynomial; the complex geometry and homotopy types are established
directly below.

The component structure of Baumslag--Solitar representation varieties has
been studied by Benyash-Krivets and Govorushko. Their 2014 note describes
the $\GL_r$-representation varieties for coprime exponents, proves that
their components have dimension $r^2$, and expresses their number in terms
of conjugacy classes of matrices with conjugate powers. Their 2016 Steklov
paper develops this description and proves rationality of the
components.\footnote{See \cite[Theorem, p.~69]{BG14} and \cite[Theorems 2.1 and 2.12]{BG16}. The target in these statements is $\GL_r$, and the standing exponent hypotheses are coprimality and the restrictions stated in those papers.}
Their separate paper on noncoprime exponents gives the component-dimension
formula
\[
 \dim H(C)=r^2+\dim Z(C^d)-\dim Z(C),
 \qquad d=\gcd(|m|,|n|),
\]
where $Z$ denotes a centralizer in $\GL_r$.\footnote{Benyash-Krivets--Govorushko \cite[Theorem, p.~53; Corollary, p.~56]{BG16nc}.}

For an affine variety $V$ of dimension $d$, let $N_i(V)$ denote
the number of its irreducible components of dimension $i$.
The \emph{profile function}, introduced by Liriano--Majewicz
\cite[Section~1, pp.~4--5]{LM}, is the ordered tuple
$P_d(V)=(N_d(V),N_{d-1}(V),\ldots,N_0(V))$.
Its entries are listed in descending order of dimension, with
$N_i(V)=0$ when no irreducible component has dimension $i$.

We build on this component-classification framework and give a direct
$\SL_2(\CC)$ calculation for the present family. This step is relevant
to the profile function: a determinant-one fiber\footnote{The fiber of a determinant map over $1$.} can have more than one
irreducible component.

Put $M=\max(|m|,|n|)$,
$\delta=\gcd(|m-n|,|m+n|)$, and $c=\gcd(2,m-n)$.
Theorem~\ref{thm:profile} gives
\begin{equation}\label{eq:profile}
 N_5(\RR)=\frac{\delta-c}{2},\qquad
 N_3(\RR)=M+3c-\delta,
 \qquad N_i(\RR)=0\quad(i\ne3,5).
\end{equation}
Thus the full profile is determined, and
\begin{equation}\label{eq:dimension}
 \dim\RR=
 \begin{cases}
 3,&\delta\le2,\\
 5,&\delta>2.
 \end{cases}
\end{equation}
In particular, dimension four does not occur in the non-Hopfian family
under consideration.

The Andreotti--Frankel--Hamm theorem places these dimensions in a
topological setting: a complex affine variety of complex dimension $D$,
with arbitrary singularities, has the homotopy type of a CW complex of
real dimension at most $D$.\footnote{Andreotti--Frankel \cite[Section 2]{AF} for the nonsingular case; Hamm \cite[Satz 1]{Hamm83} and \cite[Korollar]{Hamm86} for Stein spaces with singularities.}
Thus, when $\dim\RR=3$, ordinary integral homology can be computed from a
cellular chain complex concentrated in degrees zero through three; in
particular, $H_j(\RR;\ZZ)=0$ for $j>3$. When $\dim\RR=5$, the corresponding
CW bound is five.

The singular loci have a simpler and more uniform form. If
$Q=\{(x,y,z)\in\Aff^3_{\CC}:x^2+yz=1\}$ and $\mu=2c$, then
Theorem~\ref{thm:singular-locus} gives
\begin{equation}\label{eq:singular}
 \Singmn\cong\bigsqcup_{j=1}^{\mu}Q.
\end{equation}
Although the component profile of $\RR$ depends on $M$ and
$\delta$, the complex algebraic isomorphism type of its singular
locus depends only on the parity of $m-n$.
In particular, $\Singmn$ is smooth and pure of dimension two, whether
$\dim\RR$ is three or five. At each point of $\Singmn$, the analytic germ\footnote{An analytic germ records the complex analytic structure near
a marked point, up to shrinking the neighborhood; see Milnor
\cite{Milnor} for background.}
of $\RR$ is isomorphic to the germ of
$\CC^2\times\{(u,v)\in\CC^2:uv=0\}$ at the origin.
When $\dim\RR=5$, these singular points lie on three-dimensional components.
Since $\Singmn$ is affine of complex dimension two, the preceding CW bound gives
it a CW model of real dimension at most two. Theorem~\ref{thm:topology}
identifies its homotopy type explicitly:
$\Singmn\simeq\bigsqcup_{j=1}^{\mu}S^2$. Consequently,
\begin{equation}\label{eq:homology}
 H_j(\Singmn;\ZZ)\cong
 \begin{cases}
 \ZZ^\mu,&j=0,2,\\
 0,&\text{otherwise}.
 \end{cases}
\end{equation}
With the compactly supported Hodge--Deligne convention for the
E-polynomial,\footnote{For this convention and additivity under locally closed stratifications, see Gonz\'alez-Prieto--Mu\~noz \cite[Section 2.1]{GPM}.}
Proposition~\ref{prop:E-polynomial} yields
\begin{equation}\label{eq:E}
 E(\Singmn;u,v)=\mu\bigl((uv)^2+uv\bigr).
\end{equation}
The E-polynomial is
computed separately from an algebraic stratification of the quadrics;
compactly supported cohomology is not an invariant of ordinary homotopy
type.

To illustrate this arithmetic connection, let $R_{19}$ be the reduced
$\FF_{19}$-variety defined by $X^2Y=YX^3$ in $\SL_2^2$.
Put $k=\overline{\FF}_{19}$ and
$R_{19,k}=R_{19}\times_{\FF_{19}}k$. Then
\[
 \dim R_{19,k}=3,\qquad \dim\Sing R_{19,k}=2,
\]
and exact enumeration gives
\begin{equation}\label{eq:F19}
 \#R_{19}(\FF_{19})=34{,}200,\qquad
 \#(\Sing R_{19})(\FF_{19})=760.
\end{equation}
Here $19\nmid60$ and $6\equiv5^2\pmod{19}$, so the two singular
quadrics remain distinct and split over $\FF_{19}$.\footnote{For the
same reduced equations over $\FF_5$, one still has
$\dim R_{5,\overline{\FF}_5}=3$ and
$\dim\Sing R_{5,\overline{\FF}_5}=2$; the two trace values coalesce
at zero, so the reduced singular locus is a single smooth affine quadric surface.}
For $B(2,3)$, the complex singular surfaces have traces
$\tr Y=\pm5/\sqrt6$; choosing $\sqrt6\mapsto5$ reduces these to
$\pm1$. Thus
\[
 \Sing R_{19,k}
 =\{(I,Y):Y\in\SL_2(k),\ (\tr Y)^2=1\},
\]
the disjoint union of the trace-$1$ and trace-$-1$ quadrics.
For every $q=19^h$, $h\ge1$, their point count is $2(q^2+q)$.
More generally, the singular-locus model over
$\mathcal A=\ZZ[\sqrt6,1/30]$ is a disjoint union of two split
quadrics and has this same counting polynomial for every finite-field
specialization of $\mathcal A$. Katz's theorem therefore gives
$E(S_{2,3};u,v)=2((uv)^2+uv)$, as computed independently in
Section~\ref{sec:topology}.

\paragraph{Organization of the paper.}
Section 2 determines the dimensions and profile functions by analyzing
the scalar, semisimple, and nonsemisimple possibilities for the image of
the base generator. Section 3 identifies the singular loci and establishes
their local analytic structure. Section 4 computes their E-polynomials,
homotopy types, and integral homology. A computational supplement supplies exact enumeration code for the
$\FF_{19}$ example and records the verification procedure.

\paragraph{Computational supplement.}
The accompanying package contains a standard-library Python enumeration
over $\FF_{19}$, its output checked against the displayed counts,
and a description of the exact arithmetic and verification procedure.
The code and recorded output are supplied as ancillary files with the
arXiv submission.

\section{Irreducible components, dimensions, and profile functions}
\label{sec:components}

Throughout this section, $B(m,n)$ is non-Hopfian, $G=\SL_2(\CC)$, and
$\RR$ has its reduced affine variety structure. In particular,
$m,n\ne0$ and $|m|\ne|n|$. We retain the notation
$M=\max(|m|,|n|)$, $\delta=\gcd(|m-n|,|m+n|)$, and
$c=\gcd(2,m-n)$.

\begin{paritylemma}
The integer $\delta$ is odd if and only if $m-n$ is odd.
Consequently, $c=\gcd(2,\delta)$.
\end{paritylemma}
\begin{proof}
The integers $m+n$ and $m-n$ have the same parity. Their greatest
common divisor is therefore odd exactly when they are odd.
\end{proof}

The component descriptions of
Benyash-Krivets--Govorushko provide the framework for the calculation.\footnote{See
\cite[Theorem, p.~69]{BG14}, \cite[Theorem~2.1]{BG16}, and
\cite[Theorem, p.~53]{BG16nc}. Those statements concern $\GL_r$.
The argument below gives the $\SL_2(\CC)$ components directly, including
the two choices that arise in the nonsemisimple case.}
We first analyze the possibilities for $X$, then prove that the resulting
closed irreducible families are exactly the components.

\begin{lemma}[Spectral restriction]\label{lem:spectrum}
If $(X,Y)\in\RR$ and $\lambda,\lambda^{-1}$ are the eigenvalues of $X$,
then $\lambda^{m-n}=1$ or $\lambda^{m+n}=1$. Consequently,
$\tr X$ takes only finitely many values on $\RR$.
\end{lemma}

\begin{proof}
The defining relation gives $Y^{-1}X^mY=X^n$. Hence the unordered pairs
$\{\lambda^m,\lambda^{-m}\}$ and
$\{\lambda^n,\lambda^{-n}\}$ agree, with multiplicities.
Thus $\lambda^m$ equals either $\lambda^n$ or $\lambda^{-n}$,
which proves the alternatives. Since $m-n$ and $m+n$ are both nonzero,
the permitted eigenvalues form a finite set, as do the traces
$\lambda+\lambda^{-1}$.
\end{proof}

For $\lambda\ne\pm1$, set
$D_\lambda=\operatorname{diag}(\lambda,\lambda^{-1})$,
$T=\{\operatorname{diag}(\eta,\eta^{-1}):\eta\in\CC^\times\}$, and
$w=\left(\begin{smallmatrix}0&-1\\1&0\end{smallmatrix}\right)$.
Put $\mathcal Q_\lambda=\{Z\in G:\tr Z=\lambda+\lambda^{-1}\}$
and $F_\lambda=\{(X,Y)\in\RR:X\in\mathcal Q_\lambda\}$.

\begin{proposition}[Noncentral semisimple families]\label{prop:semisimple}
For every $\lambda\ne\pm1$ satisfying Lemma~\ref{lem:spectrum}, the
variety $F_\lambda$ is closed, smooth, and irreducible. The fiber over
$X=D_\lambda$ and the dimension of $F_\lambda$ are as follows:
\par\smallskip
{\centering\small\normalfont
\setlength{\tabcolsep}{6pt}
\renewcommand{\arraystretch}{1.03}
\begin{tabular}{@{}lcc@{}}
\toprule
Condition & Admissible $Y$ & $\dim F_\lambda$ \\
\midrule
$\lambda^{m-n}=1$, $\lambda^{m+n}\ne1$ & $T$ & $3$ \\
$\lambda^{m+n}=1$, $\lambda^{m-n}\ne1$ & $wT$ & $3$ \\
$\lambda^{m-n}=\lambda^{m+n}=1$ & $G$ & $5$ \\
\bottomrule
\end{tabular}\par}
\smallskip
Moreover, $F_\lambda=F_{\lambda^{-1}}$, and these are the only
identifications among the parameters.
\end{proposition}

\begin{proof}
Put $\tau_\lambda=\lambda+\lambda^{-1}$. Since $\tau_\lambda\ne\pm2$,
every matrix in $G$ of trace $\tau_\lambda$ is diagonalizable with
eigenvalues $\lambda,\lambda^{-1}$. Conjugacy in $\GL_2$ may be
realized in $G$ by multiplying a conjugating matrix by a scalar of
suitable determinant. Thus $\mathcal Q_\lambda$ is the conjugacy
class of $D_\lambda$, and $F_\lambda$ is a closed trace fiber in $\RR$.

Writing $X=\frac{\tau_\lambda}{2}I+
\left(\begin{smallmatrix}\alpha&\beta\\\gamma&-\alpha\end{smallmatrix}\right)$
identifies $\mathcal Q_\lambda$ with
$\alpha^2+\beta\gamma=\kappa_\lambda^2$, where
$\kappa_\lambda=(\lambda-\lambda^{-1})/2\ne0$.
Dividing the three coordinates by $\kappa_\lambda$ gives
$\mathcal Q_\lambda\cong Q$. This surface is smooth, since the
gradient $(2x,z,y)$ does not vanish on $Q$, and irreducible, since
it is the image of the irreducible group $G$ under
$g\mapsto gD_\lambda g^{-1}$.

If only $\lambda^{m-n}=1$ holds, then
$\lambda^{2m}=\lambda^{m-n}\lambda^{m+n}=\lambda^{m+n}\ne1$.
Hence $D_\lambda^m=D_\lambda^n$ has distinct eigenvalues and
centralizer $T$ in $G$.
If only $\lambda^{m+n}=1$ holds, then
$\lambda^{2m}=\lambda^{m-n}\ne1$ and
$D_\lambda^n=D_\lambda^{-m}$. Thus admissible $Y$ exchange the
two eigenlines, giving $wT$.
If both equalities hold, then $\lambda^{2m}=1$ and
$D_\lambda^m=D_\lambda^n=\lambda^mI$, so every $Y\in G$ is allowed.
For $X=gD_\lambda g^{-1}\in\mathcal Q_\lambda$, conjugation
$Y\mapsto gYg^{-1}$ identifies the fiber over $D_\lambda$ with
the fiber over $X$. The corresponding fibers are $gTg^{-1}$,
$g(wT)g^{-1}$, and $G$, respectively.

For completeness, these descriptions also establish smoothness and
irreducibility of the full families. Over $\mathcal Q_\lambda$, the
projector $P(X)=(X-\lambda^{-1}I)/(\lambda-\lambda^{-1})$ and its
complement define algebraic eigenline bundles $L_\lambda$ and
$L_{\lambda^{-1}}$. In the first case, every admissible $Y$ has the
unique form $\eta P(X)+\eta^{-1}(I-P(X))$, with $\eta\in\CC^\times$;
hence $F_\lambda\cong\mathcal Q_\lambda\times\CC^\times$.
In the second case, $Y$ interchanges the eigenlines. Choosing its
nonzero restriction $f:L_\lambda\to L_{\lambda^{-1}}$ determines
the other restriction as $-f^{-1}$, since $\det Y=1$.
Thus $F_\lambda$ is the complement of the zero section in the line
bundle $\Hom(L_\lambda,L_{\lambda^{-1}})$.
It is smooth and irreducible of dimension three. In the third case,
$F_\lambda=\mathcal Q_\lambda\times G$, of dimension five.
Finally, $\lambda+\lambda^{-1}=\lambda'+(\lambda')^{-1}$ holds
exactly when $\lambda'=\lambda$ or $\lambda^{-1}$, proving the
last assertion.
\end{proof}

Set $E=\{\varepsilon\in\{1,-1\}:\varepsilon^{m-n}=1\}$, so that
$|E|=c$. The scalar family
$A_\varepsilon=\{(\varepsilon I,Y):Y\in G\}$ will be considered
together with the nonsemisimple matrices of trace $2\varepsilon$.

\begin{proposition}[Scalar and nonsemisimple families]\label{prop:nonsemisimple}
For each $\varepsilon\in E$, the variety $A_\varepsilon$ is closed,
smooth, and irreducible of dimension three. Fix $s\in\CC^\times$
with $s^2=m/n$. For $r\in\{s,-s\}$, define
$U_{\varepsilon,r}=\{(X,Y)\in\RR:N\ne0,\ N^2=0,\ YN=rN\}$,
where $N=\varepsilon X-I$, and let
$D_{\varepsilon,r}=\overline{U_{\varepsilon,r}}$ in $\RR$.
The locus where $X$ is nonsemisimple with eigenvalue $\varepsilon$
is the disjoint union of these two irreducible three-dimensional
families. Their closures are smooth, distinct, and disjoint.
Each $D_{\varepsilon,r}$ is the total space of a line bundle over
the smooth affine quadric $\mathcal Q_r$, and
\begin{equation}\label{eq:boundary}
D_{\varepsilon,r}\setminus U_{\varepsilon,r}
=D_{\varepsilon,r}\cap A_\varepsilon
=\{\varepsilon I\}\times\mathcal Q_r.
\end{equation}
\end{proposition}

\begin{proof}
For $X=\varepsilon I$, the defining relation is equivalent to
$\varepsilon^{m-n}=1$, and then $Y$ is arbitrary. Thus
$A_\varepsilon\cong G$, proving the first assertion.

For a nonsemisimple $X$ with eigenvalue $\varepsilon$, write
$X=\varepsilon(I+N)$, where $N\ne0$ and $N^2=0$.
For $\varepsilon=\pm1$, we have
$\varepsilon^{m+n}=\varepsilon^{m-n}$, since $\varepsilon^{2n}=1$.
Lemma~\ref{lem:spectrum} therefore gives $\varepsilon^{m-n}=1$,
so $\varepsilon\in E$ and $\varepsilon^m=\varepsilon^n$.
Using $(I+N)^j=I+jN$ for every integer $j$ and cancelling this
common scalar, the relation $X^mY=YX^n$ is equivalent to $mNY=nYN$.
After conjugating $N$ to
$\left(\begin{smallmatrix}0&1\\0&0\end{smallmatrix}\right)$,
write $Y=\left(\begin{smallmatrix}\alpha&\beta\\\gamma&\omega\end{smallmatrix}\right)$,
with $\alpha,\beta,\gamma,\omega\in\CC$.
The relation gives $\gamma=0$ and $m\,\omega=n\,\alpha$, while
$\det Y=1$ gives $\alpha\,\omega=1$.
Hence $\alpha^2=m/n$, so $\alpha=\pm s$.
Taking $r=\alpha$, we obtain
$Y=\left(\begin{smallmatrix}r&\beta\\0&r^{-1}\end{smallmatrix}\right)$.
The intrinsic condition $YN=rN$ distinguishes the two families
$U_{\varepsilon,r}$ defined in the statement.

Since $|m|\ne|n|$, we have $r^2\ne1,-1$, so
$r\ne r^{-1}$ and $r+r^{-1}\ne0$.
The smooth irreducible quadric $\mathcal Q_r$ has algebraic
eigenline bundles $L_r$ and $L_{r^{-1}}$.
Consider $\Hom(L_{r^{-1}},L_r)$ over $\mathcal Q_r$.
A vector in its fiber is a map $N:L_{r^{-1}}\to L_r$, extended
by zero on $L_r$. It satisfies $N^2=0$, $YN=rN$, and $NY=r^{-1}N$.
Since $mr^{-1}=nr$, the pair $(\varepsilon(I+N),Y)$ belongs to $\RR$.

Conversely, let $(X,Y)\in U_{\varepsilon,r}$.
The preceding calculation gives $Y\in\mathcal Q_r$.
Since $mNY=nYN$ and $YN=rN$, we obtain
$NY=(nr/m)N=r^{-1}N$.
For $v\in L_r$, comparison with $Yv=rv$ gives
$(r-r^{-1})Nv=0$; hence $N|_{L_r}=0$.
Moreover, $YN=rN$ gives $N(L_{r^{-1}})\subseteq L_r$.
Thus $N$ is precisely a vector of $\Hom(L_{r^{-1}},L_r)$,
extended by zero on $L_r$.

Allowing $N=0$, the total space is the closed subvariety given by
$Y\in\mathcal Q_r$, $X=\varepsilon(I+N)$, $YN=rN$, and $NY=r^{-1}N$.
The same eigenline argument applies to every point of this closed
subvariety. The correspondence and its inverse are regular:
$N=\varepsilon X-I$, and the restrictions to the eigenlines are
expressed using the regular eigenprojectors of $Y$.
Its nonzero locus is exactly $U_{\varepsilon,r}$.
A line bundle over a smooth irreducible surface is smooth and
irreducible of dimension three, and its nonzero locus is dense.
Its total space is therefore $D_{\varepsilon,r}$, and its zero
section gives \eqref{eq:boundary}. Finally, the traces of $Y$ on
$D_{\varepsilon,s}$ and $D_{\varepsilon,-s}$ are the distinct
constants $s+s^{-1}$ and $-(s+s^{-1})$, so the closures are disjoint.
\end{proof}

\Needspace{7\baselineskip}
\begin{lemma}[Exhaustion and maximality]\label{lem:components}
The varieties $F_\lambda$, with one representative for each admissible
inverse pair $\{\lambda,\lambda^{-1}\}$ with $\lambda\ne\pm1$,
together with the $A_\varepsilon$ and $D_{\varepsilon,\pm s}$ for
$\varepsilon\in E$, are precisely the irreducible components of $\RR$.
\end{lemma}

\begin{proof}
Every matrix $X\in G$ is scalar, noncentral semisimple, or
nonsemisimple. Lemma~\ref{lem:spectrum} and
Propositions~\ref{prop:semisimple}--\ref{prop:nonsemisimple} therefore
give a finite cover of $\RR$ by the stated closed irreducible
subvarieties. It remains to exclude containments.

Distinct $F_\lambda$ lie in distinct trace fibers, and their traces
are different from $2\varepsilon$, so they are disjoint from every
scalar and nonsemisimple family. The families for different
$\varepsilon$ also lie in distinct trace fibers.
For a fixed $\varepsilon$, the two $D_{\varepsilon,\pm s}$ are
disjoint by Proposition~\ref{prop:nonsemisimple}.
Neither is contained in $A_\varepsilon$, since each contains points
with $N\ne0$. Conversely, $A_\varepsilon$ is not contained in either
closure: $\tr Y$ is unrestricted on $A_\varepsilon$ and constant
on that closure. Thus the members of this finite closed irreducible
cover are pairwise noncontained. They are exactly the maximal
irreducible closed subsets, hence exactly the irreducible components.
\end{proof}

\begin{theorem}[Component counts and profile function]\label{thm:profile}
The irreducible components of $\RR$ have dimensions three or five,
and their numbers are
\[
N_5(\RR)=\frac{\delta-c}{2},\qquad
N_3(\RR)=M+3c-\delta,\qquad
N_i(\RR)=0\quad(i\ne3,5).
\]
In particular, the profile function is completely determined by
$M$ and $\delta$.
\end{theorem}

\begin{proof}
Let $\Lambda_-=\{\lambda\in\CC^\times:\lambda^{m-n}=1\}$ and
$\Lambda_+=\{\lambda\in\CC^\times:\lambda^{m+n}=1\}$.
Their cardinalities are $|m-n|$ and $|m+n|$, and their intersection
is the group of $\delta$th roots of unity. The elements $\pm1$
that occur in either set occur in both, and are precisely the
$c$ elements of $E$.

The five-dimensional components correspond to
$(\Lambda_-\cap\Lambda_+)\setminus E$, modulo inversion.
Inversion has no fixed point on this set, so there are
$(\delta-c)/2$ such components. The three-dimensional semisimple
components correspond to the symmetric difference
$\Lambda_-\mathbin{\triangle}\Lambda_+$, again modulo inversion.
This set also has no fixed point under inversion. Its contribution is
$(|m-n|+|m+n|-2\delta)/2=M-\delta$,
using $|m-n|+|m+n|=2\max(|m|,|n|)$.
There are also $c$ scalar components and $2c$ nonsemisimple closures,
all of dimension three. Lemma~\ref{lem:components} proves that
these contributions are distinct and exhaustive, giving
$N_3(\RR)=M-\delta+c+2c$ as required.
The parity lemma gives $c=\gcd(2,\delta)$, so $M$ and $\delta$
determine the profile.
\end{proof}

\begin{corollary}[Dimension]\label{cor:dimension}
We have $\dim\RR=3$ when $\delta\le2$, and $\dim\RR=5$ when
$\delta>2$. No irreducible component has dimension four.
\end{corollary}

\begin{proof}
By the parity lemma, $\delta\le2$ is equivalent to $\delta=c$,
which is equivalent
to $N_5(\RR)=0$. Otherwise $N_5(\RR)>0$.
There is always a scalar component $A_1$ of dimension three.
The assertions follow from Theorem~\ref{thm:profile}.
\end{proof}

Listing component counts from the largest dimension down to zero,
the profiles are
$P_3(\RR)=(M+3c-\delta,0,0,0)$ when $\delta\le2$, and
$P_5(\RR)=((\delta-c)/2,0,M+3c-\delta,0,0,0)$ otherwise.
The following examples illustrate both possible dimensions:
\par\smallskip
{\centering\small
\setlength{\tabcolsep}{5pt}
\renewcommand{\arraystretch}{0.95}
\begin{tabular}{@{}lcccc@{}}
\toprule
Group & $\delta$ & $N_5$ & $N_3$ & Dimension \\
\midrule
$B(2,3)$  & 1 & 0 & 5  & 3 \\
$B(3,5)$  & 2 & 0 & 9  & 3 \\
$B(4,6)$  & 2 & 0 & 10 & 3 \\
$B(6,9)$  & 3 & 1 & 9  & 5 \\
$B(6,10)$ & 4 & 1 & 12 & 5 \\
\bottomrule
\end{tabular}\par}
\smallskip

\section{The singular locus and its local analytic structure}
\label{sec:singular-locus}

Retain the notation of Section~\ref{sec:components}, including the
fixed square root $s$ of $m/n$. For $\varepsilon\in E$ and
$r\in\{s,-s\}$, equation~\eqref{eq:boundary} gives
$A_\varepsilon\cap D_{\varepsilon,r}=\{\varepsilon I\}\times\mathcal Q_r$.
We first describe these intersections, then compute the analytic germ
of $\RR$ along them.

\begin{lemma}[Component intersections]\label{lem:intersections}
The only nonempty intersections of distinct irreducible components of
$\RR$ are the surfaces $\{\varepsilon I\}\times\mathcal Q_r$.
These $2c$ surfaces are pairwise disjoint, and each is isomorphic to
$Q=\{(x,y,z)\in\Aff^3_{\CC}:x^2+yz=1\}$.
Every point outside their union is nonsingular in $\RR$.
\end{lemma}

\begin{proof}
The components are those listed in Lemma~\ref{lem:components}.
Each $F_\lambda$ lies in a trace fiber for $X$ distinct from all the
others and from $\tr X=\pm2$; hence it meets no other component.
The families indexed by different $\varepsilon$ also have different
traces of $X$. For a fixed $\varepsilon$, Proposition~\ref{prop:nonsemisimple}
shows that $D_{\varepsilon,s}\cap D_{\varepsilon,-s}=\varnothing$ and
that $A_\varepsilon\cap D_{\varepsilon,r}=\{\varepsilon I\}\times\mathcal Q_r$.
This proves the assertion about intersections. Since
$s+s^{-1}\ne0$, the traces $s+s^{-1}$ and $-(s+s^{-1})$ are distinct,
so these surfaces are pairwise disjoint.

Projection onto the second factor identifies each surface with
$\mathcal Q_r$; the trace-quadric calculation in
Proposition~\ref{prop:semisimple} gives the asserted smoothness,
irreducibility, and isomorphism with $Q$.

Every irreducible component of $\RR$ is smooth by
Propositions~\ref{prop:semisimple} and~\ref{prop:nonsemisimple}.
A point outside the stated intersections belongs to exactly one component.
Removing the other finitely many closed components gives an open
neighborhood contained in that smooth component; the reduced structures
agree. Hence the point is nonsingular in $\RR$.
\end{proof}

\Needspace{7\baselineskip}
\begin{proposition}[Local analytic model]\label{prop:local-node}
For every $x_0\in\{\varepsilon I\}\times\mathcal Q_r$, there is an
isomorphism of complex analytic germs
\begin{equation}\label{eq:local-node}
(\RR,x_0)_{\mathrm{an}}\cong
\bigl(\CC^2\times\{(u,v)\in\CC^2:uv=0\},0\bigr).
\end{equation}
Under this isomorphism, $A_\varepsilon$ and $D_{\varepsilon,r}$
correspond to $v=0$ and $u=0$, respectively, and their intersection
corresponds to $u=v=0$.
\end{proposition}

\begin{proof}
Write $x_0=(\varepsilon I,Y_0)$. For $g\in G$, simultaneous
conjugation $(X,Y)\mapsto(gXg^{-1},gYg^{-1})$ is an algebraic
automorphism of $\RR$, fixing the scalar first coordinate
$\varepsilon I$ and preserving $A_\varepsilon$ and $D_{\varepsilon,r}$.
Indeed, the defining relation and the condition $YN=rN$ are
preserved by conjugation, so $U_{\varepsilon,r}$ and its closure
are preserved. Since $r\ne r^{-1}$, we may therefore assume
$Y_0=\operatorname{diag}(r,r^{-1})$.
For $z=(z_1,z_2)$ near $0$, set
$h(z)=\left(\begin{smallmatrix}1+z_1z_2&z_1\\z_2&1\end{smallmatrix}\right)\in G$.
For the map
$F(z_1,z_2,\zeta)=h(z)\operatorname{diag}(\zeta,\zeta^{-1})h(z)^{-1}$,
the derivatives at $(0,0,r)$ in the $z_1,z_2,\zeta$ directions are
\[
 (r^{-1}-r)E_{12},\qquad
 (r-r^{-1})E_{21},\qquad
 \operatorname{diag}(1,-r^{-2}),
\]
respectively, where $E_{ij}$ is the standard matrix unit.
These form a basis of $T_{D_r}G$, since $r\ne r^{-1}$.
The analytic inverse function theorem therefore supplies local
coordinates $(z_1,z_2,\zeta)$ on a neighborhood $U$ of $Y_0$ in $G$,
with $\zeta$ near $r$.

By Lemma~\ref{lem:spectrum}, $\tr X$ has finite image on $\RR$.
Its fiber $\tr X=2\varepsilon$ is consequently open and closed
in $\RR$, so it contains a neighborhood of $x_0$.
On this fiber put $N=\varepsilon X-I$.
The trace and determinant conditions give $\tr N=0$ and $N^2=0$.
Since $\varepsilon^m=\varepsilon^n$, the identity
$(I+N)^j=I+jN$ for every integer $j$ reduces the defining relation
to $mNY=nYN$.

Write $\widetilde N=h(z)^{-1}Nh(z)=
\left(\begin{smallmatrix}\alpha&\beta\\\gamma&-\alpha\end{smallmatrix}\right)$
and $D_\zeta=\operatorname{diag}(\zeta,\zeta^{-1})$.
The smooth ambient space is $U\times\mathfrak{sl}_2(\CC)$:
equivalently, the germ of $\CC^6$ at $(0,0,0,0,0,r)$ with
coordinates $(\alpha,\beta,\gamma,z_1,z_2,\zeta)$.
The trace fiber $\{X\in G:\tr X=2\varepsilon\}$ is the cone
$\alpha^2+\beta\gamma=0$, singular at $\varepsilon I$;
its equation is therefore included in the defining ideal.
Put $u=\zeta-r$ and work in
$\mathcal O=\CC\{z_1,z_2,u,\alpha,\beta,\gamma\}$, with $\zeta=r+u$.
The determinant-one condition and matrix relation give the ideal
\[
 I=\bigl(\alpha^2+\beta\gamma,\ (m-n)\alpha,\
         (m-n\zeta^2)\beta,\ (m\zeta^2-n)\gamma\bigr)
 \subset\mathcal O.
\]
The elements $m-n$ and $m\zeta^2-n$ are units, the latter having
nonzero value $(m^2-n^2)/n$ at $u=0$.
Moreover, $m-n\zeta^2=-n(\zeta+r)u$ and $-n(\zeta+r)$ is a unit.
Consequently, as ideals in $\mathcal O$,
\[
 I=(\alpha,\gamma,u\beta)
   =(\alpha,\gamma,u)\cap(\alpha,\gamma,\beta).
\]
The two ideals on the right are prime, so $I$ is radical.
After eliminating $\alpha,\gamma$, the remaining equation is
$u\beta=0$ in the smooth four-dimensional germ with coordinates
$(z_1,z_2,u,\beta)$.

Conversely, every sufficiently small $(z_1,z_2,\zeta,\beta)$
satisfying $(\zeta-r)\beta=0$ gives a point of $\RR$ by setting
$Y=h(z)D_\zeta h(z)^{-1}$ and
$X=\varepsilon\bigl(I+h(z)
\left(\begin{smallmatrix}0&\beta\\0&0\end{smallmatrix}\right)
h(z)^{-1}\bigr)$.
These constructions are analytic inverses. Since $I$ is radical,
it cuts out the reduced germ in the open-and-closed trace fiber.
With $v=\beta$, we obtain
\[
 \mathcal O_{\RR,x_0}^{\mathrm{an}}
 \cong\mathcal O/I
 \cong\CC\{z_1,z_2,u,v\}/(uv).
\]
This proves the claimed isomorphism of analytic germs.
The branch $v=0$ has $X=\varepsilon I$ and is $A_\varepsilon$.
On $u=0$, the points with $v\ne0$ satisfy $YN=rN$ and form
$U_{\varepsilon,r}$ locally; their closure is $D_{\varepsilon,r}$.
Their intersection is $u=v=0$, proving all assertions.
\end{proof}

\Needspace{7\baselineskip}
\begin{theorem}[The reduced singular locus]\label{thm:singular-locus}
The reduced singular locus of $\RR$ is, with
$\varepsilon\in\{\pm1\}$ and $Y\in G$,
\begin{equation}\label{eq:singular-decomposition}
\Singmn
 =\left\{(\varepsilon I,Y):
   \varepsilon^m=\varepsilon^n,\ (\tr Y)^2=\frac{(m+n)^2}{mn}\right\}
 =\bigsqcup_{\substack{\varepsilon\in E\\r\in\{s,-s\}}}
   \bigl(\{\varepsilon I\}\times\mathcal Q_r\bigr)
 \cong\bigsqcup_{j=1}^{2c}Q.
\end{equation}
It is smooth and pure of dimension two, with $\mu=2c$ irreducible
components: two when $m,n$ have opposite parity and four when they have
the same parity. Every singular point lies on exactly two
three-dimensional components of $\RR$ and has the local analytic model
\eqref{eq:local-node}. No five-dimensional component contains a singular
point of $\RR$.
\end{theorem}

\begin{proof}
Lemma~\ref{lem:intersections} shows that no point outside the displayed
union is singular. At a point of $\{\varepsilon I\}\times\mathcal Q_r$,
Proposition~\ref{prop:local-node} gives the local analytic ring
$\CC\{z_1,z_2,u,v\}/(uv)$.
It has dimension three, while its maximal ideal modulo its square
has dimension four, so the point is singular.
Algebraic and analytic nonsingularity agree, as recalled in the introduction.
This proves the equality of the reduced singular locus with the union.
For $r^2=m/n$,
$(r+r^{-1})^2=m/n+2+n/m=(m+n)^2/(mn)$.
The two choices of $r$ give the two distinct roots of this trace
equation, since $m+n\ne0$; and $\varepsilon^{m-n}=1$ is equivalent
to $\varepsilon^m=\varepsilon^n$. This proves the explicit description.

Each $\{\varepsilon I\}\times\mathcal Q_r$ is a smooth irreducible quadric surface, and
the surfaces are pairwise disjoint by Lemma~\ref{lem:intersections}.
They are closed, and thus also open in their finite union, which gives
the asserted disjoint-union isomorphism.
There are $2|E|=2c$ of them, with $c=1$ for opposite parity and $c=2$
for equal parity. The component incidence follows from the same lemma
and Proposition~\ref{prop:local-node}: only $A_\varepsilon$ and
$D_{\varepsilon,r}$ pass through a point of $\{\varepsilon I\}\times\mathcal Q_r$,
and both have dimension three. The five-dimensional components are
among the $F_\lambda$, which are smooth and meet no other component.
\end{proof}

\section{E-polynomials, homotopy types, and integral homology}
\label{sec:topology}

By Theorem~\ref{thm:singular-locus}, it suffices to study the smooth affine
quadric $Q=\{(x,y,z)\in\Aff^3_{\CC}:x^2+yz=1\}$.
We first compute its E-polynomial from an algebraic stratification, then
give an explicit description of its topology.

For a complex algebraic variety $V$, we use the compactly supported
Hodge--Deligne polynomial
\[
 E(V;u,v)=\sum_{k,i,j}(-1)^k
 h^{i,j}\bigl(H_c^k(V;\CC)\bigr)u^iv^j,
\]
where $h^{i,j}$ denotes a mixed Hodge number. We need only its additivity
under finite locally closed stratifications, its multiplicativity for
products, and the values $E(\Aff^1_{\CC};u,v)=uv$ and
$E(\CC^\times;u,v)=uv-1$; see \cite[Section~2.1]{GPM}.
Write $q=uv$ in these polynomial calculations.

\begin{proposition}[E-polynomial]\label{prop:E-polynomial}
The quadric $Q$ has $E(Q;u,v)=q^2+q$. Consequently, with $\mu=2c$,
\[
 E(\Singmn;u,v)=\mu(q^2+q).
\]
\end{proposition}

\begin{proof}
On the open subset $U=Q\cap\{y\ne0\}$, the equation determines
$z=(1-x^2)/y$. Thus
$(x,y)\mapsto(x,y,(1-x^2)/y)$ is an algebraic isomorphism
$\Aff^1_{\CC}\times\CC^\times\to U$.
The complement is the disjoint union of the two affine lines
$L_+=\{(1,0,z):z\in\CC\}$ and
$L_-=\{(-1,0,z):z\in\CC\}$.
Additivity and multiplicativity therefore give
\[
 E(Q;u,v)=E(U;u,v)+E(L_+;u,v)+E(L_-;u,v)
        =q(q-1)+2q=q^2+q.
\]
Theorem~\ref{thm:singular-locus} expresses $\Singmn$ as the disjoint
union of $\mu$ copies of $Q$, so another application of additivity
proves the second assertion.
\end{proof}

For a finite field of odd cardinality $q$, the same stratification
gives $\#Q(\FF_q)=q(q-1)+2q=q^2+q$.
Over $\mathcal A=\ZZ[\sqrt{mn},1/(2mn(m-n)(m+n))]$, the
singular-locus model is a disjoint union of $\mu$ such quadrics.\footnote{For the natural reduction over $\FF_q$, $q=p^h$ with
$p\nmid2mn(m-n)(m+n)$, the singular count is $\mu(q^2+q)$ when
$mn\in(\FF_q^\times)^2$, and zero otherwise. The square root is
included in the base $\mathcal A$ used for Katz's theorem.}
Its uniform counting polynomial is therefore $C(T)=\mu(T^2+T)$,
giving the arithmetic interpretation of this E-polynomial through
Katz's theorem.

\begin{proposition}[The quadric and the sphere]\label{prop:quadric-sphere}
The complex analytic space underlying $Q$ is diffeomorphic, as a real
smooth manifold, to the total space of the tangent bundle $TS^2$.
In particular, it admits a strong deformation retraction onto a copy
of $S^2$.
\end{proposition}

\begin{proof}
The complex linear change of coordinates
$w_1=x$, $w_2=(y+z)/2$, $w_3=(y-z)/(2i)$ carries $Q$ onto
\[
 Q'=\{w\in\CC^3:w_1^2+w_2^2+w_3^2=1\}.
\]
Write $w=\boldsymbol{\alpha}+i\boldsymbol{\beta}$ with
$\boldsymbol{\alpha},\boldsymbol{\beta}\in\mathbb R^3$, using the
Euclidean norm and inner product. The equation for $Q'$ is equivalent
to $|\boldsymbol{\alpha}|^2-|\boldsymbol{\beta}|^2=1$ and
$\boldsymbol{\alpha}\cdot\boldsymbol{\beta}=0$.
In particular, $|\boldsymbol{\alpha}|=\sqrt{1+|\boldsymbol{\beta}|^2}$
is never zero. Identify the tangent bundle with
\[
 TS^2=\{(e,\xi)\in\mathbb R^3\times\mathbb R^3:
                  |e|=1,\ e\cdot\xi=0\}.
\]
Then
\[
 \Phi:TS^2\longrightarrow Q',\qquad
 \Phi(e,\xi)=\sqrt{1+|\xi|^2}\,e+i\xi
\]
has the smooth inverse
$\boldsymbol{\alpha}+i\boldsymbol{\beta}\mapsto
(\boldsymbol{\alpha}/\sqrt{1+|\boldsymbol{\beta}|^2},\boldsymbol{\beta})$.
This proves the asserted diffeomorphism. Contracting each tangent fiber
to its zero vector gives a strong deformation retraction. Explicitly,
\[
 H_t\bigl(\Phi(e,\xi)\bigr)
   =\sqrt{1+(1-t)^2|\xi|^2}\,e+i(1-t)\xi,
 \qquad 0\le t\le1.
\]
This remains in $Q'$ because
$H_t(\Phi(e,\xi))=\Phi(e,(1-t)\xi)$ and $(e,(1-t)\xi)\in TS^2$.
At $t=0$ this is the identity; at $t=1$ it has image the real unit
sphere, and every point of that sphere is fixed for all $t$.
Transporting this retraction through the linear coordinate change
proves the result for $Q$.
\end{proof}

\begin{theorem}[Topology and homology of the singular locus]
\label{thm:topology}
Let $\mu=2c$. The singular locus has precisely $\mu$ connected
components, each simply connected, and
\[
 \Singmn\simeq\bigsqcup_{j=1}^{\mu}S^2,
 \qquad
 H_j(\Singmn;\ZZ)\cong
 \begin{cases}
 \ZZ^\mu,&j=0,2,\\
 0,&\text{otherwise}.
 \end{cases}
\]
Its ordinary integral cohomology has the same groups in the same
degrees. Its compactly supported integral cohomology is
\[
 H_c^j(\Singmn;\ZZ)\cong
 \begin{cases}
 \ZZ^\mu,&j=2,4,\\
 0,&\text{otherwise}.
 \end{cases}
\]
Both its ordinary and compactly supported Euler characteristics are
$2\mu$: four for opposite parity and eight for equal parity.
\end{theorem}

\begin{proof}
The quadrics in Theorem~\ref{thm:singular-locus} are pairwise disjoint
and open and closed in $\Singmn$. Apply
Proposition~\ref{prop:quadric-sphere} on each of them.
The resulting retractions give the stated homotopy type, connectedness,
and simple connectivity of each component.
The homology of a finite disjoint union is the direct sum of the
homology of its components. Since $H_j(S^2;\ZZ)$ is $\ZZ$ for
$j=0,2$ and zero otherwise, this proves the homology formula.
The cohomology formula follows from the integral universal coefficient
theorem, since these homology groups are free abelian.

The smooth complex surface $\Singmn$ is an oriented real manifold of
dimension four. Poincar\'e duality with compact supports, applied to
each component, gives
$H_c^j(\Singmn;\ZZ)\cong H_{4-j}(\Singmn;\ZZ)$;
see \cite[Theorem~3.35]{Hatcher}. This yields the displayed compactly
supported groups. Taking alternating sums of ranks gives $2\mu$ in
both cases, also in agreement with
$E(\Singmn;1,1)=2\mu$ from Proposition~\ref{prop:E-polynomial}.
\end{proof}

The E-polynomial and the homotopy type are computed independently,
using the algebraic stratification and the deformation retraction,
respectively. The compactly supported groups
above use the manifold structure and Poincar\'e duality in addition to
the homotopy calculation.

\section*{Acknowledgements and computational background}
This paper grew out of decades of investigations into the representation
varieties of non-Hopfian Baumslag--Solitar groups and years of experimentation
with SageMath, Singular, Macaulay2, GAP, and related tools.\footnote{For the
	software, see \cite{Sage,Singular,M2,GAP}.}
The experiments were conducted on a Linux host machine, with local VMs also 
deployed for many of the calculations.

OpenAI was instrumental in preparing the manuscript, helping with
its organization and the clarity and flow of its exposition, and supporting
extensive discussions of the soundness of the results. Other containerized local 
AI agents assisted in locating relevant literature from around the web.  The author takes responsibility for the mathematical content.

Although the arguments here admit much greater generality, the present scope was
chosen to keep the exposition short, naive, and focused, without invoking
Grothendieckian formalism. A broader framework may become appropriate in
future work on representation varieties of groups in algebraic groups over
fields of characteristic $p>0$.

\par\bigskip
\noindent
\begin{minipage}{\textwidth}
\small
Foro Cultural 81\\
Centro Hist\'orico\\
Guanajuato, Gto., 36000, Mexico\\
Email: \href{mailto:SAL21458@yahoo.com}{\nolinkurl{SAL21458@yahoo.com}}
\end{minipage}
\end{document}